\documentclass{amsart}
\usepackage{amssymb}
\usepackage{eucal}
\theoremstyle{plain}
\newtheorem{theorem}{Theorem}

\newtheorem{lemma}{Lemma}

\theoremstyle{definition}

\theoremstyle{remark}

\newtheorem{remark}{Remark}
\numberwithin{equation}{section}
\newcommand{\e}{\epsilon}

\newcommand{\de}{\delta}
\newcommand{\f}{\theta^{x,t}_\delta}
\newcommand{\R}{\mathbb R}

\newcommand{\Z}{\mathbb Z}
\newcommand{\C}{\mathbb C}

\newcommand{\Rn}{\mathbb R^n}

\newcommand{\Rm}{\mathbb R^{n+1}}
\begin{document}

\title[The sharp Hardy Uncertainty Principle for Sch\"odinger evolutions]{The sharp Hardy Uncertainty Principle for Sch\"odinger evolutions}
\author{L. Escauriaza}
\address[L. Escauriaza]{UPV/EHU\\Dpto. de Matem\'aticas\\Apto. 644, 48080 Bilbao, Spain.}
\email{luis.escauriaza@ehu.es}
\thanks{The first and fourth authors are supported  by MEC grant, MTM2007-62186, the second and third authors by NSF grants DMS-0456583 and DMS-0456833 respectively}
\author{C. E. Kenig}
\address[C. E. Kenig]{Department of Mathematics\\University of Chicago\\Chicago, Il. 60637 \\USA.}
\email{cek@math.uchicago.edu}
\author{G. Ponce}
\address[G. Ponce]{Department of Mathematics\\
University of California\\
Santa Barbara, CA 93106\\
USA.}
\email{ponce@math.ucsb.edu}
\author{L. Vega}
\address[L. Vega]{UPV/EHU\\Dpto. de Matem\'aticas\\Apto. 644, 48080 Bilbao, Spain.}
\email{mtpvegol@lg.ehu.es}
\begin{abstract}
We give a new proof of Hardy's uncertainty principle, up to the end-point case, which is only based on calculus. The method allows us to extend Hardy's uncertainty principle to Schr\"odinger equations with non-constant coefficients. We also deduce optimal Gaussian decay bounds for solutions to these Schr\"odinger equations.
\end{abstract}
\maketitle
\begin{section}{Introduction}\label{S: Introduction} 
In this paper we continue the study initiated in \cite{kpv02}, \cite{ekpv06}, \cite{ekpv08a} and \cite{ekpv08b} on unique continuation properties of solutions of Schr\"odinger evolutions
\begin{equation}\label{E: 1.1}
\partial_t u=i\left(\triangle u+V(x,t) u\right)\ ,\ \text{in}\  \R^n\times [0,T].
\end{equation}

The goal is to obtain sharp and non-trivial sufficient conditions on a solution $u$, the potential $V$ and the  behavior of the solution at two different times, $t_0<t_1$, which guarantee that $u\equiv 0$ in $\Rn\times [t_0,t_1]$.

One of our motivations  comes from a well known result due to G. H. Hardy \cite[pp. 131]{StSh}, which concerns the decay of a function $f$ and its Fourier transform, 
\[\hat f(\xi)=(2\pi)^{-\frac n2}\int_{\Rn}e^{-i\xi\cdot x}f(x)\,dx.\]

\emph{If 
$f(x)=O(e^{-|x|^2/\beta^2})$, $\hat f(\xi)=O(e^{-4|\xi|^2/\alpha^2})$ and  $1/\alpha\beta>1/4$, then $f\equiv 0$. Also, if $1/\alpha\beta=1/4$, $f$ is a constant multiple of $e^{-|x|^2/\beta^2}$.}

\vspace{0,1 cm}
As far as we know, the only known proof of this result and its variants uses complex analysis (the Phragm\'en-Lindel\"of principle). There has also been considerable interest in a better understanding of this result and on extensions of it to other settings: \cite{CoPr}, \cite{Hor2}, \cite{SiSu}, \cite{bonamie} and \cite{bonamie2}.

This result can be rewritten in terms of the free solution of the Schr\"odinger equation in $\Rn\times (0,+\infty)$, $i\partial_tu+\triangle u=0$, with initial data $f$,
\begin{equation*}
u(x,t)= (4\pi it)^{-\frac n 2} \int_{\Rn}e^{\frac{ i |x-y|^2}{4t}}f(y)\, dy= \left(2\pi it\right)^{-\frac n2}e^{\frac{i|x|^2}{4t}}\widehat{e^{\frac{i|\,\cdot\,|^2}{4t}}f}\left(\frac x{2t}\right),
\end{equation*}
in the following way: 
\vspace{0,1 cm}

\emph{If $u(x,0)=O(e^{-|x|^2/\beta^2})$, $u(x,T)=O(e^{-|x|^2/\alpha^2})$ and $T/\alpha\beta> 1/4$, then $u\equiv 0$. Also, if $T/\alpha\beta=1/4$, $u$ has as initial data a constant multiple of $e^{-\left(1/\beta^2+i/4T\right)|y|^2}$.}
\vspace{0,1 cm}

The corresponding results in terms of $L^2$-norms, established in \cite{SiSu}, are the following: 

\vspace{0,1 cm}
\emph{If $e^{|x|^2/\beta^2}f$, $e^{4|\xi |^2/\alpha^2}\widehat f$ are in $L^2(\Rn)$ and $1/\alpha\beta\ge 1/ 4$, then $f\equiv 0$.}

\vspace{0,1 cm}
\emph{If $e^{|x|^2/\beta^2}u(x,0)$, $e^{|\xi |^2/\alpha^2}u(x,T)$ are in $L^2(\Rn)$ and $T/\alpha\beta\ge 1/4$, then $u\equiv 0$.}
\vspace{0,2 cm}

In \cite {ekpv08b} we proved a uniqueness result in this direction for bounded potentials $V$ verifying, $V(x,t)=V_1(x)+V_2(x,t)$, with $V_1$ real-valued and
\[\sup_{[0,T]}\|e^{T^2|x|^2/\left(\alpha t+\beta\left(T-t\right)\right)^2}V_2(t)\|_{L^\infty(\Rn)}<+\infty\]
or
\begin{equation*}
\lim_{R\rightarrow +\infty}\int_0^T\|V(t)\|_{L^\infty(\Rn\setminus B_R)}\,dt =0.
\end{equation*}

More precisely, we proved that the only solution $u$ to \eqref{E: 1.1} in $C([0,T], L^2(\Rn))$, which verifies 
\begin{equation}\label{E: condicion fundamental}
\|e^{|x|^2/\beta^2}u(0)\|_{L^2(\Rn)}+\|e^{|x|^2/\alpha^2}u(T)\|_{L^2(\Rn)}<+\infty
\end{equation}
is the zero solution, when $T/\alpha\beta>1/ 2$ and $V$ verifies one of the above conditions.

This linear result was then applied to show that  two regular solutions $u_1$ and $u_2$ of non-linear equations of the type
\begin{equation}
\label{E: NLS}
i\partial_tu+\triangle u=F(u,\overline u), \ \text{in}\ \Rn\times [0,T]
\end{equation}
and for very general non-linearities $F$, must agree in $\Rn\times [0,T]$, when $u=u_1-u_2$ satisfies \eqref{E: condicion fundamental}. This replaced the assumption that the solutions coincide on large sub-domains of $\Rn$ at two different times, which was studied in \cite{kpv02} and \cite{Ioke04}, and showed that (weaker) variants of Hardy's Theorem hold even in the context of non-linear Schr\"odinger evolutions.

The main results in this paper improve the results in \cite{ekpv06}, \cite{ekpv08b} and show that the  optimal version of Hardy's Uncertainty Principle in terms of $L^2$-norms, as established in \cite{SiSu}, holds for solutions to \eqref{E: 1.1}, when $T/\alpha\beta>1/4$ and for many general bounded potentials, while it fails for some complex-valued potentials in the end-point case, $T/\alpha\beta=1/4$. As a by product of our argument, sharp improvements of Gaussian decay estimates are also obtained.

\begin{theorem}\label{T: hardytimeindepent}
Assume that $u$ in $C([0,T]),L^2(\Rn))$ verifies
\begin{equation*}
\label{E: 1.111}\partial_tu=i\left(\triangle u+V(x,t)u\right),\ \text{in}\ \Rn\times [0,T],
\end{equation*}
$\alpha$ and $\beta$ are positive, $T/\alpha\beta > 1/4$, $\|e^{|x|^2/\beta^2}u(0)\|_{L^2(\Rn)}$ and $\|e^{|x|^2/\alpha^2}u(T)\|_{L^2(\Rn)}$ are both finite, the potential $V$ is  bounded and either, $V(x,t)=V_1(x)+V_2(x,t)$, with $V_1$ real-valued and 
\begin{equation*}
\sup_{[0,T]}\|e^{T^2|x|^2/\left(\alpha t+\beta \left(T-t\right)\right)^2}V_2(t) \|_{L^\infty(\Rn)} < +\infty
\end{equation*}
or $\lim_{R\rightarrow +\infty}\|V\|_{L^1([0,T], L^\infty(\Rn\setminus B_R)}=0$. Then, $u\equiv 0$.
\end{theorem}

\begin{theorem}\label{T: hardytimeindepent2}
Assume that $T/\alpha\beta=1/4$. Then, there is a smooth complex-valued potential $V$ verifying 
\begin{equation*}
|V(x,t)|\lesssim\frac 1{1+|x|^2},\ \text{in}\ \R^n\times [0,T]
\end{equation*}
 and a nonzero smooth function $u$ in $C^\infty([0,T],\mathcal S(\Rn))$ such that
\begin{equation*}
\partial_tu=i\left(\triangle u+V(x,t)u\right),\ \text{in}\ \Rn\times [0,T]
\end{equation*}
and $\|e^{|x|^2/\beta^2}u(0)\|_{L^2(\Rn)}$ and $\|e^{|x|^2/\alpha^2}u(T)\|_{L^2(\Rn)}$ are both finite.
\end{theorem}

Our proof of Theorem \ref{T: hardytimeindepent} does not use any complex analysis, it only uses calculus!  It provides the first proof (up to the end-point) that we know of Hardy's uncertainty principle for the Fourier transform, without the use of complex analysis.
 
 As a by product, we derive the following optimal interior estimate for the Gaussian decay of solutions to \eqref{E: 1.1}.

\begin{theorem}\label{T: lamejora}
Assume that $u$ and $V$ verify the hypothesis in Theorem \ref{T: hardytimeindepent} and $T/\alpha\beta\le 1/4$. Then,

\begin{multline*}
\sup_{[0,T]}\|e^{a(t)|x|^2}u(t)\|_{L^2(\Rn)} +\| \sqrt{t(T-t)}\nabla \left(e^{\left(a(t)+\frac{i\dot a(t)}{8a(t)}\right)|x|^2}u\right)\|_{L^2(\Rn\times [0,T])}\\
\le N\left[\|e^{|x|^2/\beta^2}u(0)\|_{L^2(\Rn)}+ \|e^{|x|^2/\alpha^2}u(T)\|_{L^2(\Rn)}\right],
\end{multline*}
where
\[a(t)=\frac {\alpha\beta RT}{2\left(\alpha t+\beta (T-t)\right)^2+2R^2\left(\alpha t - \beta (T-t)\right)^2}\ ,\]
$R$ is the smallest root of the equation 
\begin{equation*}
\frac T{\alpha\beta}=\frac R{2\left(1+R^2\right)}
\end{equation*}
and $N$ depends on $T$, $\alpha$, $\beta$ and the conditions on the potential $V$ in Theorem \ref{T: hardytimeindepent}.
\end{theorem} 

Observe that $1/a(t)$ is convex and attains its minimum value in the interior of $[0,T]$, when
\begin{equation*}
|\alpha-\beta|<R^2\left(\alpha+\beta\right).
\end{equation*}

To understand why Theorem \ref{T: lamejora} is optimal, observe that
\begin{equation}\label{E: el enemigo}
u_R(x,t)=R^{-\frac n2}\left(t-\tfrac iR\right)^{-\frac n2}e^{-\frac{|x|^2}{4i(t-\frac iR)}}=\left(Rt-i\right)^{-\frac n2}e^{-\frac{(R-iR^2t)}{4(1+R^2t^2)}\,|x|^2},
\end{equation}
is a free wave (i.e. $V\equiv 0$, in \eqref{E: 1.1}) satisfying in $\Rn\times [-1,1]$ the corresponding time translated conditions in Theorem \ref{T: lamejora} with $T=2$ and
\begin{equation*}
\frac1{\beta^2}=\frac1{\alpha^2}=\mu=\frac R{4\left(1+R^2\right)}\le\frac 18\, .
\end{equation*}
Moreover
\begin{equation*}
\frac R{4\left(1+R^2t^2\right)}\, ,
\end{equation*}
is increasing in the $R$-variable, when  $0<R\le 1$ and $-1\le t\le 1$. See also Lemma \ref{L: transformada}.

As a direct  consequence of  Theorem \ref{T: hardytimeindepent} we get the following application concerning the uniqueness of solutions for non-linear equations of the form \eqref{E: NLS}.

\begin{theorem}
\label{Theorem 1...}

Let $u_1$ and $u_2$ be strong solutions in $C([0,T],H^k(\R^n))$ of the equation (1.3)  with $k\in \Z^+$, $k>n/2$,
$F:\C^2\to \C$, $F\in C^k$  and $F(0)=\partial_uF(0)=\partial_{\bar u}F(0)=0$. If there are $\alpha$ and $\beta$ positive  with $T/\alpha \beta>1/4$ such that
\begin{equation*}
e^{|x|^2/\beta^2}\left(u_1(0)-u_2(0)\right)\ \text{and}\ e^{|x|^2/\alpha^2}\left(u_1(T)-u_2(T)\right) 
\end{equation*}
are in $L^2(\Rn)$, then $u_1\equiv u_2$.

\end{theorem}

Notice that the condition, $T/\alpha\beta>1/4$, is independent of the size of the potential or the dimension and that we do not assume any decay of the gradient, neither of the solutions or of time-independent potentials or any regularity of the potentials.  

Our improvement for the results in \cite{ekpv06} and \cite{ekpv08b} comes from  a better understanding of the solutions to \eqref{E: 1.1}, which have Gaussian decay. We started the study of this particular type of solutions in \cite{ekpv08a}, where we considered free waves. The improvement of the latter results is a consequence of the possibility of extending the following outline of a strategy to prove Theorem \ref{T: hardytimeindepent} for free waves to the non-free wave cases:

First, by a suitable change of variables based on the conformal or Appell transform (See Lemma \ref{L: transformada}), it suffices to prove Theorem \ref{T: hardytimeindepent}, when $u$ in $C([-1,1], L^2(\Rn))$ is a solution of 
\begin{equation}\label{E: free wave}
\partial_tu-i\triangle u=0,\ \text{in}\ \Rn\times [-1,1]
\end{equation}
and
\begin{equation}\label{E: decaimineto}
\|e^{\mu |x|^2}u(-1)\|_{L^2(\Rn)}+\|e^{\mu |x|^2}u(1)\|_{L^2(\Rn)}<+\infty,
\end{equation}
for some $\mu >0$. Our strategy consists of showing that either $u\equiv 0$ or there is a function $\theta_{R}: [-1,1]\longrightarrow [0,1]$ such that
\begin{equation*}
\|e^{\frac{R|x|^2}{4\left(1+R^2t^2\right)}}u(t)\|_{L^2(\Rn)}\le \|e^{\mu |x|^2}u(-1)\|_{L^2(\Rn)}^{\theta_{R}(t)}\|e^{\mu |x|^2}u(1)\|_{L^2(\Rn)}^{1-\theta_{R}(t)},
\end{equation*}
where $R$ is the smallest root of the equation
\[\mu =\frac{R}{4\left(1+R^2\right)}\ .\]
Thus, we obtain the optimal improvement of the Gaussian decay of a free wave verifying \eqref{E: decaimineto}  and we derive that $\mu \le 1/8$, when $u$ is not zero. 

The proof of these facts relies on new logarithmic convexity properties of free waves verifying \eqref{E: decaimineto} and on those already established in \cite{ekpv08b}. In \cite[Theorem 3]{ekpv08b}, the positivity of the space-time commutator of the symmetric and skew-symmetric parts of the operator,
\begin{equation*}
e^{\mu |x|^2}\left(\partial_t-i\triangle\right)e^{-\mu |x|^2},
\end{equation*}
is used to show that  $\|e^{\mu |x|^2}u(t)\|_{L^2(\Rn)}$ is logarithmically convex in $[-1,1]$. In particular, that
\begin{equation*}
\|e^{\mu |x|^2}u(t)\|_{L^2(\Rn)}\le \|e^{\mu |x|^2} u(-1)\|_{L^2(\Rn)}^{\frac{1-t}2}\|e^{\mu |x|^2} u(1)\|_{L^2(\Rn)}^{\frac{1+t}2},
\end{equation*}
when, $-1\le t\le 1$.

Beginning from this fact we set, $a_1\equiv \mu$, and we begin a constructive procedure, where at the $k$th step, we construct $k$ smooth even functions, $a_i:[-1,1]\longrightarrow (0,+\infty)$, $1\le i\le k$, such that
\begin{equation*}
\mu\equiv a_1<a_2<\dots<a_k,\ \text{in}\ (-1,1),
\end{equation*}
\begin{equation*}
F(a_i)> 0,\ \text{in}\  [-1,1],\ a_i(1)=\mu,\ i=1,\dots,k,
\end{equation*}
where
\begin{equation*}
F(a)=\frac 1a\left(\ddot a-\frac{3\dot a^2}{2a\,}+32a^3\right)
\end{equation*}
and functions $\theta_i:[-1,1]\longrightarrow [0,1]$, $1\le i\le k$, such that
\begin{equation}\label{E: algoagradable}
\|e^{a_i(t) |x|^2}u(t)\|_{L^2(\Rn)}\le \|e^{\mu |x|^2}u(-1)\|_{L^2(\Rn)}^{\theta_i(t)}\|e^{\mu |x|^2}u(1)\|_{L^2(\Rn)}^{1-\theta_i(t)},\quad -1\le t\le1.
\end{equation}
 
These estimates are proved from the construction of the functions $a_i$, while the method strongly relies on the following formal convexity properties of free waves: 
\begin{equation}\label{E: algo fundamental}
\partial_t\left(\tfrac 1a\partial_t\log{H_b}\right)\ge -\frac{2\ddot b^2|\xi|^2}{F(a)},
\end{equation}
\begin{equation*}
\partial_t\left(\tfrac 1a\partial_tH\right)\ge \epsilon_a\int_{\Rn}e^{a|x|^2}\left(|\nabla u|^2+|x|^2|u|^2\right)\,dx,
\end{equation*}
where 
\begin{equation*}
H_b(t)=\|e^{a(t)|x+ b(t)\xi|^2}u(t)\|_{L^2(\Rn)}^2\ ,\  H(t)=\|e^{a(t)|x|^2}u(t)\|_{L^2(\Rn)}^2,
\end{equation*}
$\xi\in\Rn$ and $a, b: [-1,1]\longrightarrow\R$ are smooth functions with 
\begin{equation*}
a> 0,\quad F(a)>0,\ \text{in}\ [-1,1].
\end{equation*}

Once the $k$th step is completed, we take  $a=a_k$ in \eqref{E: algo fundamental} with a certain choice of $b=b_k$, verifying $b(-1)=b(1)=0$ and then, a certain test is performed. When the answer to the test is positive, it follows that $u\equiv 0$. Otherwise, the logarithmic convexity associated to \eqref{E: algo fundamental} allows us to find a new smooth function $a_{k+1}$ in $[-1,1]$ with
\begin{equation*}
a_1<a_2<\dots<a_k<a_{k+1},\ \text{in}\ (-1,1),
\end{equation*}
and verifying the same properties as $a_1,\dots,a_k$.

When the process is infinite, we have \eqref{E: algoagradable} for all $k\ge 1$ and there are two possibilities: either $\lim_{k\to +\infty}a_k(0)=+\infty$ or $\lim_{k\to +\infty}a_k(0)<+\infty$. The first case and \eqref{E: algoagradable} implies that $u\equiv 0$, while in the second, the sequence $a_k$ is shown to converge to an even function $a$ verifying
\begin{equation*}
\begin{cases}
\ddot a-\frac{3\dot a^2}{2a\,\,}+32a^3=0,\ \text{in}\ [-1,1],\\
a(1)=\mu.
\end{cases}
\end{equation*}
Because
\begin{equation*}
\frac{R}{4\left(1+R^2t^2\right)}\, ,\quad R\in\R,
\end{equation*}
are all the possible even solutions of this equation, $a$ must be one of them and 
\[\mu =\frac{R}{4\left(1+R^2\right)}\, ,\]
for some $R>0$. In particular, $u\equiv 0$, when $\mu >1/8$.

The proof of Theorem \ref{T: hardytimeindepent} for non-zero potentials $V$ relies on extending the above convexity properties to the non-free case. The path that goes from the formal level to a rigorous one is not an easy one. In fact in \cite[\S 6]{ekpv08b}, we gave explicit examples of functions $a(t)$ such that $\log H$ is formally convex, when 
\[H(t)=\|e^{a(t)|x|^2}u(t)\|_{L^2(\Rn)}^2\] 
but for which, the corresponding inequalities  lead  to false statements: {\it all free waves verifying \eqref{E: decaimineto} for some $\mu >0$ are identically zero}. Therefore most parts of this paper, as those in \cite{ekpv08b}, are devoted to making rigorous the above formal arguments.

\end{section}
\begin{section}{A few Lemmas}\label{S:  Algunos Lemmas}
In the sequel \[\left(f,g\right)=\int_{\Rn}f\overline g\,dx\ ,\ \|f\|^2=\left(f,f\right).\]

The following formal identities or inequalities appeared or were proved within the proof of \cite[Lemma 2]{ekpv08b}.

\begin{lemma}\label{L: freq1}
$\mathcal S$ is a symmetric operator, $\mathcal A$ is skew-symmetric, both are allowed to depend on the time variable, $f(x,t)$ is a reasonable function, 
\begin{equation}\label{E: unas definiciones}
H(t)=\left( f, f\right)\ ,\ \partial_t\mathcal S=\mathcal S_t\ ,\ D(t)=\left(\mathcal Sf,f\right)\ \text{and}\ N(t)=\frac{D(t)}{H(t)}\,.
\end{equation}
Then,
\begin{equation}\label{E: primera derivada}
\dot H(t)=2\text{\it Re}\left(\partial_tf-\mathcal Sf-\mathcal Af,f\right)+2\left(\mathcal Sf,f\right), 
\end{equation}
\begin{multline}
\label{E: derivadasegunda}
\ddot H\ge 2\partial_t\text{\it Re}\left(\partial_tf-\mathcal Sf-\mathcal Af,f\right)+2\left(\mathcal S_tf+\left[\mathcal S,\mathcal A\right]f,f\right)-\|\partial_tf-\mathcal Af-\mathcal Sf\|^2,
\end{multline}
\begin{equation}\label{E: logaritmo21}
\dot D(t)\ge \left(\mathcal S_tf+\left[\mathcal S,\mathcal A\right]f,f\right)-\frac 12\|\partial_tf-\mathcal Af-\mathcal Sf\|^2
\end{equation}
and
\begin{equation}\label{E: desigualdaddefrecuenciacasiexacta}
\dot N(t)\ge \left(\mathcal S_tf +\left[\mathcal S,\mathcal A\right]f, f\right)/H- \|\partial_tf-\mathcal Af-\mathcal Sf\|^2/\left(2H\right).
\end{equation}
\end{lemma}

Lemma \ref{L: aproximaci—n de convexidad logar'tmica} shows how to find possible convexity or  log-convexity properties of $H(t)$ with respect to a new and possibly unknown variable $s$, which is related the original time variable $t$ by the ordinary differential equation
\begin{equation*}
\frac{dt}{ds}=\gamma(t).
\end{equation*}
\begin{lemma}\label{L: aproximaci—n de convexidad logar'tmica} Assume that $\mathcal S$, $\mathcal A$ and $f$ are as above, $\epsilon >0$, and $\gamma: [c,d]\longrightarrow (0,+\infty)$ and $\psi:[c,d]\longrightarrow [0,+\infty)$ are smooth functions satisfying
\begin{equation}\label{E: condicioncomutadorgeneralizado}
\left(\gamma\,\mathcal S_t f(t)+\gamma\left[\mathcal S,\mathcal A\right]f(t)+\dot\gamma\,\mathcal Sf(t),f(t)\right)\ge -\psi(t)\, H(t), \ \text{when}\ c\le t\le d.
\end{equation}
Then,
\begin{align}\label{E: una propiedad de convexidad logar'tmica}
H(t)+\epsilon\le \left(H(c)+\epsilon\right)^{\theta(t)}\left(H(d)+\epsilon\right)^{1-\theta(t)}e^{2T(t)+M_\epsilon(t)+2N_\e}, \ \text{when}\ c\le t\le d,
\end{align}
where $T$ and $M_\e$ verify
\begin{equation*}
\partial_t\left(\gamma\,\partial_tT\right)=-\psi,\ \text{in}\ [c,d],\quad T(c)=T(d)=0,
\end{equation*}

\begin{equation*}
\partial_t\left(\gamma\,\partial_tM_\e\right)=-\gamma\,\frac{\|\partial_tf-\mathcal Sf-\mathcal Af\|^2}{H+\e}\, ,\ \text{in}\ [c,d],\quad M_\e(c)=M_\e(d)=0,
\end{equation*}
\begin{equation*}
N_\e=\int_c^d\left |\text{\it Re}\ \frac{\left(\partial_sf(s)-\mathcal Sf(s)-\mathcal Af(s),f(s)\right)}{H(s)+\e}\right |\,ds
\end{equation*}
and 
\begin{equation*}\label{E: el jodidoexponente}
\theta(t)=\frac{{\int_t^d}\frac{ds}{\gamma}}{\int_c^d\frac{ds}{\gamma}}\, .
\end{equation*}
Moreover,
\begin{align}\label{E: formula para la derivada segunda generalizada}
\partial_t\left(\gamma\,\partial_t H-\gamma{\it Re}\left(\partial_tf-\mathcal Sf-\mathcal Af,f\right)\right)&+\gamma\,\|\partial_tf-\mathcal Sf-\mathcal Af\|^2\\&\ge
2\left(\gamma\,\mathcal S_tf+\gamma\left[\mathcal S,\mathcal A\right]f+\dot\gamma\,\mathcal Sf,f\right).\notag 
\end{align}
\end{lemma}
\begin{proof}
From \eqref{E: primera derivada}
\begin{align}\label{E: primera derivada1}
\partial_t\log{\left(H+\epsilon\right)}-2\text{\it Re}\,\frac{\left(\partial_tf-\mathcal Sf-\mathcal Af,f\right)}{H+\e}&=\frac{2\left(\mathcal Sf,f\right)}{H+\e}\\&=\frac{2D}{H+\e}.\notag 
\end{align}
The differentiation of the second identity in \eqref{E: primera derivada1} gives
\begin{equation*}
\partial_t^2\log{\left(H+\epsilon\right)}-2\partial_t\text{\it Re}\,\frac{\left(\partial_tf-\mathcal Sf-\mathcal Af,f\right)}{H+\e}=\frac{2\e\dot D}{\left(H+\e\right)^2}+\frac{2H^2\dot N}{\left(H+\e\right)^2}\, ,
\end{equation*}
with $D$ and $N$ as defined in \eqref{E: unas definiciones}, and from \eqref{E: logaritmo21} and \eqref{E: desigualdaddefrecuenciacasiexacta}
\begin{multline}\label{E: otraderivadasegunda}
\partial_t^2\log{\left(H+\epsilon\right)}-2\partial_t\text{\it Re}\,\frac{\left(\partial_tf-\mathcal Sf-\mathcal Af,f\right)}{H+\e}\\\ge
\frac{2\left(\gamma\,\mathcal S_t+\gamma\left[\mathcal S,\mathcal A\right]f+\dot\gamma\,\mathcal Sf,f\right)}{H+\e}-\gamma\,\frac{\|\partial_tf-\mathcal Sf-\mathcal Af\|^2}{H+\e}\,.
\end{multline}
Multiply the first identity in \eqref{E: primera derivada1} by $\dot\gamma$, \eqref{E: otraderivadasegunda} by $\gamma$ and add up the corresponding identity and inequality to obtain the inequality

\begin{multline*}
\partial_t\left(\gamma\,\partial_t\log{\left(H+\epsilon\right)}-2\gamma \text{\it Re}\,\frac{\left(\partial_tf-\mathcal Sf-\mathcal Af,f\right)}{H+\e}\right)\\\ge
\frac{2\left(\gamma\,\mathcal S_t+\gamma\left[\mathcal S,\mathcal A\right]f+\dot\gamma\,\mathcal Sf,f\right)}{H+\e}-\gamma\,\frac{\|\partial_tf-\mathcal Sf-\mathcal Af\|^2}{H+\e}\,.
\end{multline*}
This and \eqref{E: condicioncomutadorgeneralizado} show that
\begin{multline*}
\partial_t\left(\gamma\,\partial_t\log{\left(H+\epsilon\right)}-2\gamma\text{\it Re}\,\frac{ \left(\partial_tf-\mathcal Sf-\mathcal Af,f\right)}{H+\e}\right)\\\ge -2\psi-\gamma\,\frac{\|\partial_tf-\mathcal Sf-\mathcal Af\|^2}{H+\e}\,.
\end{multline*}
Thus,
\begin{equation*}\label{E: algo de monotonia}
\partial_t\left(\gamma\left(\partial_t\log{\left(H+\epsilon\right)}-2\partial_tT-\partial_tM_\e-2\text{\it Re}\,\frac{\left(\partial_tf-\mathcal Sf-\mathcal Af,f\right)}{H+\e}\right)\right)\ge 0,
\end{equation*}
in $[c,d]$. The monotonicity associated to  
this inequality shows that

\begin{align*}
\frac 1{\gamma(\tau)}[\partial_s \log{\left (H(s)+\e\right)}-2\partial_sT(s)-\partial_s& M_\e(s)] \\
-\frac 2{\gamma(\tau)}\,\text{\it Re}\,&\frac{\left(\partial_sf(s)-\mathcal Sf(s)-\mathcal Af(s),f(s)\right)}{H(s)+\e}
\end{align*}
\begin{align*}
\le\frac 1{\gamma(s)}[ \partial_\tau \log{\left(H(\tau)+\e\right)}-2\partial_\tau T(\tau)-\partial_\tau& M_\e (\tau)]\\
-\frac 2{\gamma(s)}\,\text{\it Re}\,&\frac{\left(\partial_\tau f(\tau)-\mathcal Sf(\tau)-\mathcal Af(\tau),f(\tau)\right)}{H(\tau)+\e}\, ,
\end{align*}
when $c\le s \le \tau\le d$, and the integration of this inequality for $(s,\tau)$ in $[c,t]\times [t,d]$ and the boundary conditions satisfied by $T$ and $M_\e$, imply that the \lq\lq logarithmic convexity\rq\rq\  type inequality \eqref{E: una propiedad de convexidad logar'tmica} holds, when $c\le t\le d$.

To derive \eqref{E: formula para la derivada segunda generalizada} multiply the identity \eqref{E: primera derivada} by $\dot\gamma$, the inequality \eqref{E: derivadasegunda} by $\gamma$ and add up the corresponding identity and inequality.
\end{proof}
A calculation (See also formulae (2.12), (2.13) and (2.14) in \cite{ekpv08b}) shows that given smooth functions $a:[c,d]\longrightarrow (0,+\infty)$ and $b: [c,d]\longrightarrow\Rn$,  
\begin{equation*}
e^{a(t)|x+ b(t)|^2}\left(\partial_t-i\triangle\right)e^{-a(t)|x+ b(t)|^2}=\partial_t-\mathcal S-\mathcal A,
\end{equation*}
where $\mathcal S$ and $\mathcal A$ are respectively symmetric and skew-symmetric operators on $\Rn$, given by the formulae, 
\begin{align}\label{E: parte simetrica y antisimetrica}
&\mathcal S=-2i\left(2a\left(x+b\right)\cdot\nabla+an\right)+\dot a|x+b|^2+2a\dot b\cdot\left(x+b\right),\\
&\mathcal A=i\left(\triangle +4a^2|x+b|^2\right).\notag
\end{align}
Moreover,
\begin{multline}\label{E: el conmutador}
\mathcal S_t+\left[\mathcal S,\mathcal A\right]=-8a\triangle -2i\left(\left(4\dot a\left(x+b\right)+4a\dot b\right)\cdot\nabla+2\dot an\right)\\+\left(\ddot a+32a^3\right)|x+b|^2+\left(4\dot a\dot b+2a\ddot b\right)\cdot\left(x+b\right)+2a|\dot b|^2.
\end{multline}

In Lemma \ref{L: el calculo del conmutador mas largo} we calculate a lower bound for the self-adjoint operator \[\gamma\,\mathcal S_t+\gamma\left[\mathcal S,\mathcal A\right]+\dot\gamma\,\mathcal S,\]
 when $\gamma: [c,d]\longrightarrow (0,+\infty)$ is a smooth function.

\begin{lemma}\label{L: el calculo del conmutador mas largo} Let $a$,\ $\gamma$ and $b$ be as above and assume that
\begin{equation}\label{E: formula para el termino del cuadrado}
F(a,\gamma)=\gamma\left(\ddot a+32a^3-\frac{3\dot a^2}{2a}-\frac {a}2\left(\frac{\dot a}a+\frac{\dot\gamma}\gamma\right)^2\right)>0,\ \text{in}\ [c,d].
\end{equation}
Then, if $\mathcal I$ denotes the identity operator,
\begin{equation*}
\gamma\,\mathcal S_t+\gamma\left[\mathcal S,\mathcal A\right]+\dot\gamma\,\mathcal S\ge -\frac{\gamma^2a^2|\ddot b|^2}{F(a,\gamma)}\,\mathcal I,
\end{equation*}
for each time $t$ in $[c,d]$.
\end{lemma}

\begin{proof}
 From \eqref {E: parte simetrica y antisimetrica}
, \eqref{E: el conmutador} and the identity
\begin{equation*}
\left(\gamma\,\mathcal S_tf+\gamma\left[\mathcal S,\mathcal A\right]f+\dot\gamma\,\mathcal Sf,f\right)=\text{\it Re}\left(\gamma\,\mathcal S_tf+\gamma\left[\mathcal S,\mathcal A\right]f+\dot\gamma\,\mathcal Sf,f\right),
\end{equation*}
we have
\begin{multline*}
\left(\gamma\,\mathcal S_tf+\gamma\left[\mathcal S,\mathcal A\right]f+\dot\gamma\,\mathcal Sf,f\right)=\int_{\Rn}\left(32\gamma a^3+\gamma\ddot a+\dot\gamma\dot a\right)|x+b|^2|f|^2\,dx\\
+\int_{\Rn}\left[\left(4\gamma\dot a\dot b+2\gamma a\ddot b +2a\dot\gamma\dot b\right)\cdot \left(x+b\right)+2\gamma a|\dot b|^2\right]|f|^2\,dx\\
+\int_{\Rn}8\gamma a|-i\nabla f|^2+2\text{\it Re}\left(-i\nabla f\right)\cdot\overline{\left(4\gamma a \dot b f\right)}\\+\int_{\Rn}2\text{\it Re}\left(-i\nabla f\right)\cdot\overline{\left(\left(2\dot\gamma a+4\gamma\dot a\right)\left(x+b\right) f\right)}\, dx,
\end{multline*}
when $f$ is in $S(\Rn)$. Completing the square corresponding to the first and second terms in the third line above, we get
\begin{multline*}
\left(\gamma\,\mathcal S_tf+\gamma\left[\mathcal S,\mathcal A\right]f+\dot\gamma\,\mathcal Sf,f\right)=\int_{\Rn}\left(32\gamma a^3+\gamma\ddot a+\dot\gamma\dot a\right)|x+b|^2|f|^2\,dx\\
+\int_{\Rn}\left(4\gamma\dot a\dot b+2\gamma a\ddot b +2a\dot\gamma\dot b\right)\cdot \left(x+b\right)|f|^2\,dx\\
+\int_{\Rn}8\gamma a|-i\nabla f+\tfrac{\dot b}2\,f|^2+2\text{\it Re}\left(-i\nabla f\right)\cdot\overline{\left(\left(2\dot\gamma a+4\gamma\dot a\right)\left(x+b\right) f\right)}\, dx.
\end{multline*}

Rewriting $-i\nabla f$ in the second term in the third line above as 
\[\left(-i\nabla f+\tfrac{\dot b}2\,f\right)- \tfrac{\dot b}2\,f,\]
gives the formula
\begin{multline*}
\left(\gamma\,\mathcal S_tf+\gamma\left[\mathcal S,\mathcal A\right]f+\dot\gamma\,\mathcal Sf,f\right)=\int_{\Rn}\left(32\gamma a^3+\gamma\ddot a+\dot\gamma\dot a\right)|x+b|^2|f|^2\,dx\\
+\int_{\Rn}2\gamma a\ddot b\cdot \left(x+b\right)|f|^2\,dx\\
+\int_{\Rn}8\gamma a|-i\nabla f+\tfrac{\dot b}2\,f|^2+2\text{\it Re}\left(-i\nabla f+\tfrac{\dot b}2\,f\right)\cdot\overline{\left(\left(2\dot\gamma a+4\gamma\dot a\right)\left(x+b\right) f\right)}\, dx.
\end{multline*}

Next, we complete the square corresponding to the two terms in the third line above and find that
\begin{multline}\label{E: casi la formula}
\left(\gamma\,\mathcal S_tf+\gamma\left[\mathcal S,\mathcal A\right]f+\dot\gamma\,\mathcal Sf,f\right)=\int_{\Rn}8\gamma a|-i\nabla f+\tfrac{\dot b}2\,f+\left(\tfrac{\dot a}{2a}+\tfrac{\dot\gamma}{4\gamma}\right)\left(x+b\right)f|^2\, dx\\+\int_{\Rn}\left[F(a,\gamma)|x+b|^2+2\gamma a\ddot b\cdot \left(x+b\right)\right]|f|^2\,dx,
\end{multline}
where
\begin{align*}
F(a,\gamma)&=\ddot a\gamma+32\gamma a^3+\dot\gamma\dot a-\frac {\gamma a}2\left(\frac{2\dot a}a+\frac{\dot\gamma}\gamma\right)^2\\
&=\gamma\left(\ddot a+32a^3-\frac{3\dot a^2}{2a\,}-\frac {a}2\left(\frac{\dot a}a+\frac{\dot\gamma}\gamma\right)^2\right).
\end{align*}

Finally, we complete the square corresponding to the terms in the second line of \eqref{E: casi la formula} to obtain that
\begin{multline}\label{E: vaya formulita}
\left(\gamma\,\mathcal S_tf+\gamma\left[\mathcal S,\mathcal A\right]f+\dot\gamma\,\mathcal Sf,f\right)=\int_{\Rn}8\gamma a|-i\nabla f+\tfrac{\dot b}2\,f+\left(\tfrac{\dot a}{2a}+\tfrac{\dot\gamma}{4\gamma}\right)\left(x+b\right)f|^2\, dx\\+\int_{\Rn}F(a,\gamma)|x+b+\tfrac{a\gamma\ddot b}{F(a,\gamma)}|^2|f|^2\,dx-\frac{\gamma^2a^2|\ddot b|^2}{F(a,\gamma)}\int_{\Rn}|f|^2,
\end{multline}
when $f$ is in $\mathcal S(\Rn)$ and $c\le t\le d$, which proves Lemma \ref{L: el calculo del conmutador mas largo}.
\end{proof}
In the sequel we set 
\begin{equation*}
F(a)=F(a,\tfrac 1a)= \frac 1a\left(\ddot a+32a^3-\frac{3\dot a^2}{2a\,}\right).
\end{equation*}

The main Lemma used here in justifying the formal calculations is the following.
\begin{lemma}\label{L: justificaci—n de las cuentas} Assume that $u$ in $C([-1,1],L^2(\Rn))$ verifies
\begin{equation}\label{E: la ecuaci—n de u}
\partial_t u=i\left(\triangle u+V(x,t) u\right)\ ,\ \text{in}\  \R^n\times [-1,1]\ ,\ \sup_{[-1,1]}\|e^{\mu|x|^2}u(t)\|<+\infty,
\end{equation}
for some $\mu>0$ and some bounded complex-valued potential $V$. Also assume that $a:[-1,1]\longrightarrow (0,+\infty)$ and $b: [-1,1]\longrightarrow\Rn$ are smooth, $b(-1)=b(1)=0$, $a$ is even, $\dot a\le 0$ in $[0,1]$, $a(1)=\mu$, $a\ge \mu$ and $F(a)>0$ in $[-1,1]$, and that
\begin{equation}\label{E: una condicion fundamental}
\sup_{[-1,1]}\|e^{\left(a(t)-\e\right)|x|^2}u(t)\|<+\infty,
\end{equation}
for all $\e>0$. Then,
\begin{equation}\label{E: convexiad}
\|e^{a(t)|x+b(t)|^2}u(t)\|\le e^{T(t)+2\|V\|_\infty+\frac 14\|V\|^2_\infty}\sup_{[-1,1]}\|e^{\mu|x|^2}u(t)\|,\ -1\le t\le 1
\end{equation}
where $T$ verifies
\begin{equation}\label{E: laformuladeT}
\begin{cases}
\partial_t\left(\tfrac 1a\,\partial_tT\right)=-\frac{|\ddot b|^2}{F(a)},\ \text{in}\ [-1,1],\\
T(-1)=T(1)=0.
\end{cases}
\end{equation}

Moreover, there is $\eta_a>0$, such that
\begin{align}\label{E: mas regularidad2}
&\|\sqrt{1-t^2}\,\,\nabla(e^{(a+\frac{\dot a i}{8a})|x|^2}u)\|_{L^2(\Rn\times [-1,1])}\\
&+\eta_a\|\sqrt{1-t^2}\,e^{a(t)|x|^2}\nabla u\|_{L^2(\Rn\times [-1,1])}\le e^{2\|V\|_\infty+\frac 14\|V\|^2_\infty}\sup_{[-1,1]}\|e^{\mu|x|^2}u(t)\|.\notag
\end{align}

\end{lemma}

\begin{proof} 
Extend $u$ to $\Rm$ as $u\equiv 0$, when $|t|>1$. For $\e>0$, set 
\begin{equation*}
a_\e(t)=a(t)-\e,\ g_\e(x,t)=e^{a_\e(t)|x|^2}u(x,t),\ f_\e(x,t)=e^{a_\e(t)|x+b(t)|^2}u(x,t).
\end{equation*}
Then, $f_\e$ is in $L^\infty([-1,1],L^2(\Rn))$ and verifies 
\begin{equation}\label{E: la ecuacion de f}
\partial_tf_\e-\mathcal S_\e f_\e-\mathcal A_\e f_\e=iV(x,t)f_\e,
\end{equation}
in the sense of distributions in $\Rn\times (-1,1)$, with $\mathcal S_\e$, $\mathcal A_e$ defined by formulae \eqref{E: parte simetrica y antisimetrica} but with $a$ replaced by $a_\e$:
\begin{equation}\label{E: soluciondebil}
\int f_\e\overline{\left(-\partial_s\xi-\mathcal S_\e\xi+\mathcal A_\e\xi\right)}\,dyds=i\int Vf_\e\overline\xi\,dyds,\ \text{for all}\ \xi\in C^\infty_0(\Rn\times (-1,1)).
\end{equation}

Let $\theta$ in $C^\infty(\Rm)$ be a standard mollifier supported in the unit ball of $\Rm$ and for $0<\de\le \frac 14$ set, $g_{\e,\delta}=g_\e\ast\theta_\delta$, $f_{\e,\delta}=f_\e\ast\theta_\delta$ and
\begin{equation*}
\theta^{x,t}_\de(y,s)=\delta^{-n-1}\theta(\tfrac{x-y}\delta\, ,\tfrac{t-s}\delta).
 \end{equation*}
 Then, $f_{\e,\delta}$ is in $C^\infty([-1,1],\mathcal S(\Rn))$ and
\begin{align}\label{E: formulalarga}
&\left(\partial_tf_{\e,\de}-\mathcal S_\e f_{\e,\de}-\mathcal A_\e f_{\e,\de}\right)(x,t)\\
&=\int f_\e\overline{\left(-\partial_s\f-\mathcal S_\e\f+\mathcal A_\e\f\right)}dyds\notag\\
&+\int f_\e\left[\left(\dot a_\e(s)+4ia_\e^2(s)\right)|y+b(s)|^2-\left(\dot a_\e(t)+4ia_\e^2(t)\right)|x+b(t)|^2\right]\f dyds\notag\\
&+\int f_\e\left[2a_\e(s)\dot b(s)\cdot \left(y+b(s)\right)-2a_\e(t)\dot b(t)\cdot \left(x+b(t)\right)\right]\f dyds\notag\\
&+\int f_\e\left[4i\left(a_\e(s)\left(y+b(s)\right)-a_\e(t)\left(x+b(t)\right)\right)\cdot\nabla_y\f\right]dyds\notag\\
&+\int 2in f_\e\left(a_\e(s)+a_\e(t)\right)\f\,dyds,\notag
\end{align}
when $x$ is in $\Rn$ and $-1+\de\le t\le 1-\de$. The last identity and \eqref{E: soluciondebil} give,
\begin{equation}\label{E: controldelerror}
\left(\partial_tf_{\e,\de}-\mathcal S_\e f_{\e,\de}-\mathcal A_\e f_{\e,\de}\right)(x,t)=i\left(Vf_\e\right)\ast\theta_\de(x,t)+A_{\e,\de}(x,t)+B_{\e,\de}(x,t),
\end{equation}
in $\Rn\times [-1+\delta,1-\delta]$, where $A_{\e,\de}$ and $B_{\e,\de}$ denote respectively the sum of the second and third integrals and of the fourth and fifth in \eqref{E: formulalarga}. Moreover, there is $N_{a,b,\e}$ such that

\begin{align*}
&|A_{\e,\de}(x,t)|\le \delta^{-n} N_{a,b,\e}\int_{t-\de}^{t+\de}\int_{B_\delta(x)}|e^{\left(a(s)-\frac \e 2\right)|y|^2}|u(y,s)|\,dyds,\\
&|B_{\e,\de}(x,t)|\le \delta^{-1-n} N_{a,b,\e}\int_{t-\de}^{t+\de}\int_{B_\delta(x)}|e^{\left(a(s)-\frac \e 2\right)|y|^2}|u(y,s)|\,dyds,
\end{align*}
when $(x,t)$ is in $\Rn\times [-1+\de,1-\de]$, which implies that
\begin{align}
\sup_{[-1+\de,1-\de]}&\|A_{\e,\de}(t)\|_{L^2(\Rn\times [-1+\delta,1-\de])}\le \de N_{a,b,\e}\sup_{[-1,1]}\|e^{\left(a(t)-\frac\e 2\right)|x|^2}u(t)\|,\label{E: controldeA}\\
\|B_{\e,\de}&\|_{L^2(\Rn\times [-1+\delta,1-\de])}\le N_{a,b,\e}\sup_{[-1,1]}\|e^{\left(a(t)-\frac\e2\right)|x|^2}u(t)\|.\label{E:control malo de B}
\end{align}
We also have,
\begin{equation}\label{E: el termino digno}
\sup_{[-1,1]}\|\left(Vf_\e\right)\ast\theta_\de(t)\|\le \|V\|_\infty\sup_{[-1,1]}\|e^{\left(a(t)-\tfrac\e2\right)|x|^2}u(t)\|.
\end{equation}

Clearly, $g_{\e,\delta}$ also verifies \eqref{E: controldelerror} with the corresponding $A_{\e,\de}$ and $B_{\e,\de}$ verifying  \eqref{E: controldeA} and \eqref{E:control malo de B}. Just set $b\equiv 0$ in the definitions of $\mathcal S_\e$, $\mathcal A_\e$ and replace $f_\e$ by $g_\e$ in $A_{\e,\de}$ and $B_{\e,\de}$. We can also replace $f_\e$ by $g_\e$ in \eqref{E: el termino digno}.

From the hypothesis on $a$ and $b$, there is $\e_a>0$ such that 
\begin{equation}\label{E: otrapropiedadagradable}
F(a_\e)\ge \tfrac12\, F(a),\ \text{in}\  [-1,1],\ \text{when}\ 0<\e\le\e_a
\end{equation}
and for such an  $\e>0$, it is possible to find $\delta_\e>0$, with $\delta_\e$ approaching zero as $\e$ tends to zero, such that
\begin{equation}\label{E: algunas propiedades}
\left(a(t)-\tfrac\e2\right)|x|^2,\   \left(a(t)-\tfrac \e2\right)|x+b(t)|^2\le \mu |x|^2,\ \text{when}\ x\in\Rn, \ 1-\delta_\e\le |t|\le 1.
\end{equation}

	From now on, always assume that $0<\e\le\e_a$ and $0<\de\le\de_\e$. Then, the calculations leading to the identity \eqref{E: primera derivada1} are justified, when $[c,d]=[-1+\de,1-\de]$, $\gamma=\frac 1{a_\e}$, $\mathcal S=\mathcal S_\e$, $\mathcal A=\mathcal A_\e$, $H_{\e,\de}(t)=\|g_{\e,\de}(t)\|^2$ and we get
\begin{align}\label{E: conelmalditoepsi–on}
&\partial_t\left(\tfrac 1{a_\e}\,\partial_t H_{\e,\de}-\tfrac 1{a_\e}{\it Re}\left(\partial_t\,g_{\e,\de}-\mathcal S_\e\, g_{\e,\de}-\mathcal A_\e\, g_{\e,\de}\,,\,g_{\e,\de}\right)\right) \\ &+\tfrac 1{a_\e}\,\|\partial_t\,g_{\e,\de}-\mathcal S_\e\, g_{\e,\de}-\mathcal A_\e\, g_{\e,\de}\|^2\notag\\ &\ge
2\left(\tfrac 1{a_\e}\,\mathcal S_{\e t}\,g_{\e,\de}+\tfrac 1{a_\e}\left[\mathcal S_\e,\mathcal A_{\e}\right]g_{\e,\de}-\tfrac{\dot a_\e}{a^2_\e}\,\mathcal S_\e\, g_{\e,\de}\, ,\,g_{\e,\de}\right). \notag
\end{align}
Moreover, the identity \eqref{E: vaya formulita} in Lemma \ref{L: el calculo del conmutador mas largo} with $\gamma=\frac 1{a_\e}$, $b\equiv 0$ and \eqref{E: otrapropiedadagradable} show that
\begin{multline}\label{E: vaya formulita2}
\left(\tfrac 1{a_\e}\,\mathcal S_{\e t}f+\tfrac 1{a_\e}\left[\mathcal S_{\e},\mathcal A_\e\right]f-\tfrac{\dot a_\e}{a^2_\e}\,\mathcal S_\e f,f\right)\\\ge\int_{\Rn}8|-i\nabla f+\tfrac{\dot a_\e}{4a_\e}\,xf|^2+F(a_\e)|x|^2|f|^2\,dx\\
\ge \int_{\Rn} |\nabla(e^{\frac{\dot a_\e i}{8a_\e}|x|^2}f)|\,dx+\eta_a^2\int_{\Rn}|\nabla f|^2+|x|^2|f|^2\,dx,
\end{multline}
when $f\in\mathcal S(\Rn)$. The multiplication of the inequality \eqref{E: conelmalditoepsi–on} by $\left(1-\delta_\e\right)^2-t^2$, integration by parts, \eqref{E: controldeA}, \eqref{E:control malo de B}, \eqref{E: el termino digno} and the fact that $V$ is bounded in $L^\infty$ imply that
\begin{equation}\label{E: mas regularidad}
\|\sqrt{\left(1-\delta_\e\right)^2-t^2}\,\nabla g_{\e,\de}\|_{L^2(\Rn\times [-1+\de_\e,1-\de_\e])}\le N_{a,\e}\, .
\end{equation}
Analogously,
\begin{equation*}
\|\sqrt{\left(1-\delta_\e\right)^2-t^2}\,\nabla f_{\e,\de}\|_{L^2(\Rn\times [-1+\de_\e,1-\de_\e])}\le N_{a,b,\e}
\end{equation*}
and after letting $\delta$ tend to zero, we find that
\begin{equation}\label{E: controlgradiente3}
\|\sqrt{\left(1-\delta_\e\right)^2-t^2}\,\nabla f_{\e}\|_{L^2(\Rn\times [-1+\de_\e,1-\de_\e])}\le N_{a,b,\e}\, .
\end{equation}
The latter makes it possible to write the error term $B_{\e,\de}$ as
\begin{equation*}
\int \left[4i\left(a_\e(t)\left(x+b(t)\right)-a_\e(s)\left(y+b(s)\right)\right)\cdot\nabla_yf_\e+2in\left(a_\e(t)-a_\e(s)\right)f_\e\right]\f dyds
\end{equation*}
and derive that
\begin{equation}\label{E: algomasiendoacero}
\|B_{\e,\de}\|_{L^2(\Rn\times [-1+\de_\e,1-\de_\e])}\le \delta N_{a,b,\e}\, ,\ \text{when} \ 0<\de\le\de_\e.
\end{equation}

Recalling Lemma \ref{L: el calculo del conmutador mas largo},  apply now the estimate \eqref{E: una propiedad de convexidad logar'tmica} on logarithmic convexity to $f_{\e,\de}$, with $H_{\e,\delta}(t)=\|f_{\e,\de}(t)\|^2$,  $[c,d]=[-1+\de_\e,1-\de_\e]$, $\gamma=\frac 1{a_\e}$, $\mathcal S=\mathcal S_\e$ and $\mathcal A=\mathcal A_\e$, and from \eqref{E: otrapropiedadagradable} and \eqref{E: algunas propiedades}, we get
\begin{equation}\label{E: casialli}
H_{\e,\de}(t)\le \left(\sup_{[-1,1]}\|e^{\mu|x|^2}u(t)\|+\e\right)^2e^{2T_\e(t)+M_{\e,\de}(t)+2N_{\e,\de}}, \ \text{when}\ |t|\le 1-\de_\e,
\end{equation}
where $T_\e$ and $M_{\e,\de}$ verify
\begin{equation*}
\begin{cases}
\partial_t\left(\tfrac 1{a_\e}\,\partial_tT_\e\right)=-\frac{|\ddot b|^2}{F(a_\e)},\ \text{in}\ [-1+\delta_\e,1-\delta_\e],\\\ T_\e(-1+\delta_\e)=T_\e(1-\delta_\e)=0,
\end{cases}
\end{equation*}

\begin{equation*}
\begin{cases}
\partial_t\left(\tfrac 1{a_\e}\,\partial_tM_{\e,\de}\right)=-\tfrac 1{a_\e}\,\frac{\|\partial_tf_{\e,\de}-\mathcal S_\e f_{\e,\de}-\mathcal A_\e f_{\e,\de}\|^2}{H_{\e,\de}+\e}\, ,\ \text{in}\ [-1+\de_\e,1-\delta_\e],\\
M_{\e,\delta}(-1+\de_\e)=M_{\e,\de}(1-\de_\e)=0,
\end{cases}
\end{equation*}
and
\begin{equation*}
N_{\e,\de}=\int_{-1+\delta_\e}^{1+\delta_\e}\frac{\|\partial_sf_{\e,\de}(s)-\mathcal S_\e f_{\e,\de}(s)-\mathcal A_\e f_{\e,\de}(s)\|}{\sqrt{H_{\e,\de}(s)+\e}}\,ds.
\end{equation*}

The equation \eqref{E: la ecuacion de f} verified by $f_\e$, \eqref{E: una condicion fundamental}, \eqref{E: controlgradiente3} and the formulae \eqref{E: parte simetrica y antisimetrica} show that $f_\e$ is in $C((-1,1),L^2(\Rn))$ and $H_{\e,\de}$ converges uniformly on compact sets of $(-1,1)$ to $H_\e(t)=\|f_\e(t)\|^2$. From \eqref{E: controldeA}, \eqref{E: algomasiendoacero} and letting first $\delta$ tend to zero and then $\e$ tend to zero in \eqref{E: casialli}, we get
\begin{align*}
\|e^{a(t)|x+b(t)|^2}u(t)\|^2\le \sup_{[-1,1]}\|e^{\mu|x|^2}u(t)\|^2e^{2T(t)+M(t)+4\|V\|_\infty}, \ \text{when}\ |t|\le 1,
\end{align*} 
where $T$ was defined in \eqref{E: laformuladeT} and
\begin{equation*}
\begin{cases}
\partial_t\left(\tfrac 1{a}\,\partial_tM\right)=-\tfrac 1{a}\,\|V\|_\infty^2,\\
M(-1)=M(1)=0.
\end{cases}
\end{equation*}

Because $M$ is even,
\begin{equation*}
M(t)=\|V\|_\infty^2\int_t^1\int_0^s\frac{a(s)}{a(\tau)}\,d\tau ds,\ \text{in}\ [0,1],  
\end{equation*}
and the monotonicity of $a$ in $[0,1]$ implies that $M(t)\le \|V\|_\infty^2/2$, in $[-1,1]$, which proves \eqref{E: convexiad}. The inequality
\begin{equation}\label{E: final}
\|e^{a(t)|x|^2}u(t)\|\le e^{2\|V\|_\infty+\frac 14\|V\|^2_\infty}\sup_{[-1,1]}\|e^{\mu|x|^2}u(t)\|,\ -1\le t\le 1
\end{equation} 
and the reasoning leading to \eqref{E: mas regularidad} can be repeated again, but now using  \eqref{E: final}, together with \eqref{E: controldeA}, \eqref{E: algomasiendoacero} and \eqref{E: vaya formulita2}, to show that
\begin{align*}
&\|\sqrt{\left(1-\delta_\e\right)^2-t^2}\,\nabla(e^{\frac{\dot a_\e i}{8a_\e}|x|^2}g_{\e,\de})\|_{L^2(\Rn\times [-1+\de_\e,1-\de_\e])}\\
&+\eta_a\|\sqrt{\left(1-\delta_\e\right)^2-t^2}\,\nabla g_{\e,\de}\|_{L^2(\Rn\times [-1+\de_\e,1-\de_\e])}\\&+\eta_a\|\sqrt{\left(1-\delta_\e\right)^2-t^2}\, x\, g_{\e,\de}\|_{L^2(\Rn\times [-1+\de_\e,1-\de_\e])}\\&\le e^{2\|V\|_\infty+\frac 14\|V\|^2_\infty}\sup_{[-1,1]}\|e^{\mu|x|^2}u(t)\|+\delta N_{a,\e}.
\end{align*}

Letting first $\de$ tend to zero, and then $\e$ tend to zero, we get \eqref{E: mas regularidad2}, which proves Lemma \ref{L: justificaci—n de las cuentas}.

\end{proof}
\end{section}
\begin{section}{Proofs of Theorems \ref{T: hardytimeindepent}, \ref{T: hardytimeindepent2} and \ref{T: lamejora}}\label{Proofs of Theorems}
We first recall the following Lemma proved in \cite[Lemma 5]{ekpv08b}. It is useful to reduce the case of different Gaussian decays at two distinct times to the the case of the same Gaussian decay.
\begin{lemma}\label{L: transformada} Assume that  $\alpha$ and $\beta$ are positive, $\gamma\in\R$ and that $u$ in $C([0,1], L^2(\Rn))$ verifies
\begin{equation*}
\partial_su=i\left(\triangle_y u+V(y,s)u\right)\ ,\ \text{in}\  \R^n\times [0,1].
\end{equation*}
Set
\[\widetilde u(x,t)=\left(\tfrac{\sqrt{\alpha\beta}}{\alpha(1-t)+\beta t}\right)^{\frac n2}u\left(\tfrac{\sqrt{\alpha\beta}\, x}{\alpha(1-t)+\beta t}, \tfrac{\beta t}{\alpha(1-t)+\beta t}\right)e^{\left(\alpha-\beta\right) |x|^2/4i(\alpha(1-t)+\beta t)}.\]
Then, $\widetilde u$ is in $C([0,1], L^2(\Rn))$ and verifies
 \begin{equation*}
\partial_t\widetilde u=i\left(\triangle_x \widetilde u+\widetilde V(x,t)\widetilde u\right)\ ,\ \text{in}\  \R^n\times [0,1],
\end{equation*}
with
\begin{equation*}
\widetilde V(x,t)=\tfrac{\alpha\beta}{\left(\alpha(1-t)+\beta t\right)^2}\,V\left(\tfrac{\sqrt{\alpha\beta}\, x}{\alpha(1-t)+\beta t}, \tfrac{\beta t}{\alpha(1-t)+\beta t}\right).
\end{equation*}
Moreover,
\begin{equation}\label{E: que maravilladecambio}
\|e^{\gamma |x|^2}\widetilde u(t)\| = \|e^{\gamma\alpha\beta |y|^2/(\alpha s+\beta(1-s))^2}u(s)\|,
\end{equation}
when $s=\tfrac{\beta t}{\alpha(1-t)+\beta t}$.
\end{lemma}
\begin{proof}[Proof of Theorems \ref{T: hardytimeindepent} and \ref{T: lamejora}] Let $u$ satisfy the conditions in Theorem \ref{T: hardytimeindepent} and set $u_T(x,t)=T^{\frac n4}u(\sqrt Tx,Tt)$. We have,
\begin{equation*}
\partial_tu_T=i\left(\triangle u_T+V_T\left(x,t\right)\right),\ \text{in}\ \Rn\times [0,1],
\end{equation*}
with $V_T(x,t)=TV(\sqrt Tx,Tt)$ and
\begin{equation*}
\|e^{|x|^2/\beta^2}u(0)\|+\|e^{|x|^2/\alpha^2}u(T)\|=\|e^{T|x|^2/\beta^2}u_T(0)\|+\|e^{T|x|^2/\alpha^2}u_T(1)\|.
\end{equation*}
From \cite[Theorem 3]{ekpv08b}, when $V$ verifies the first condition or \cite[Theorem 5]{ekpv08b}, when $V$ verifies the second condition in Theorem \ref{T: hardytimeindepent}, we know that
\begin{equation*}
\sup_{[0,T]}\|e^{T^2|x|^2/\left(\alpha t+\beta\left(T-t\right)\right)^2}u(t)\|=\sup_{[0,1]}\|e^{T|x|^2/\left(\alpha t+\beta\left(1-t\right)\right)^2}u_T(t)\|<+\infty.
\end{equation*}
In fact, there it is shown that
\begin{multline}\label{E: la primera acotaci—n}
\sup_{[0,1]}\|e^{T|x|^2/\left(\alpha t+\beta\left(1-t\right)\right)^2}u_T(t)\|\\\le N\left[\|e^{T|x|^2/\beta^2}u_T(0)\|+\|e^{T|x|^2/\alpha^2}u_T(1)\|\right],
\end{multline}
where $N$ depends on $\alpha$, $\beta$ and the conditions imposed on the potential $V$ in Theorem \ref{T: hardytimeindepent}. From Lemma \ref{L: transformada} 
\begin{equation*}
\widetilde u(x,t)=\left(\tfrac{\sqrt{\alpha\beta}}{\alpha(1-t)+\beta t}\right)^{\frac n2}u_T\left(\tfrac{\sqrt{\alpha\beta}\, x}{\alpha(1-t)+\beta t}\, , \tfrac{\beta t}{\alpha(1-t)+\beta t}\right)e^{\left(\alpha-\beta\right) |x|^2/4i\left(\alpha(1-t)+\beta t\right)},
\end{equation*}
verifies 
\begin{equation*}
\partial_t\widetilde u=i(\triangle \widetilde u+\widetilde V_T(x,t)\widetilde u)\ ,\ \text{in}\  \R^n\times [0,1],
\end{equation*}
with 
\begin{equation*}
\widetilde V_T(x,t)=\tfrac{\alpha\beta T}{\left(\alpha(1-t)+\beta t\right)^2}\,V\left(\tfrac{\sqrt{\alpha\beta T}\, x}{\alpha(1-t)+\beta  t}, \tfrac{\beta Tt}{\alpha(1-t)+\beta t}\right),
\end{equation*}
and from \eqref{E: que maravilladecambio},
\begin{equation*}
\sup_{[0,1]}\|e^{T|x|^2/\alpha\beta}\widetilde u(t)\|=\sup_{[0,1]}\|e^{T|x|^2/\left(\alpha t+\beta\left(1-t\right)\right)^2}u_T(t)\|<+\infty.
\end{equation*}
Finally, $v(x,t)=2^{-\frac n4}\widetilde u(\frac x{\sqrt 2}, \frac{1+t}2)$, verifies 
\begin{equation*}
\partial_tv=i\left(\triangle v+\mathcal V(x,t)v\right),\ \text{in}\ \Rn\times [-1,1],
\end{equation*}
with $\mathcal V(x,t)=\frac 12\widetilde V_T(\frac x{\sqrt 2},\frac{1+t}2)$ and 
\begin{equation}\label{E: primer paso}
\sup_{[-1,1]}\|e^{\mu|x|^2}v(t)\|=\sup_{[0,1]}\|e^{T|x|^2/\alpha\beta}\,\widetilde u(t)\|<+\infty,\ \text{with} \ \mu=\tfrac T{2\alpha\beta}\, .
\end{equation}

Abusing notation, we replace in what follows $v$ and $\mathcal V$ by $u$ and $V$, and set, $a_1(t)\equiv \mu$. We then begin an inductive procedure, where at the $k$th step we have constructed $k$ smooth even functions, $a_j:[-1,1]\longrightarrow (0,+\infty)$ and numbers $\eta_{a_j}>0$, such that
\begin{equation*}
a_1<a_2<\dots<a_k,\ \text{in}\ (-1,1),
\end{equation*}
\begin{equation}\label{E: algunas propiedades mas}
\dot a_j\le 0\ \text{in}\ [0,1],\ F(a_j)> 0,\ \text{in}\  [-1,1],\ a_j(1)=\mu,
\end{equation}
\begin{equation}\label{E: algoagradable222}
\sup_{[-1,1]}\|e^{a_j(t) |x|^2}u(t)\|\le e^{2\|V\|_\infty+\frac 14\|V\|_\infty^2}\sup_{[-1,1]}\|e^{\mu |x|^2}u(t)\|,
\end{equation}
and
\begin{multline}\label{E: algoagradable2222}
\|\sqrt{1-t^2}\,\,\nabla(e^{(a_j+\frac{i\dot a_j}{8a_j})|x|^2}u)\|_{L^2(\Rn\times [-1,1])}\\
+\eta_{a_j}\|\sqrt{1-t^2}\,e^{a_j(t)|x|^2}\nabla u\|_{L^2(\Rn\times [-1,1])}\\\le e^{2\|V\|_\infty+\frac 14\|V\|_\infty^2}\sup_{[-1,1]}\|e^{\mu|x|^2}u(t)\|,
\end{multline}
when $j=1,\dots,k$.
The case $k=1$ follows from \eqref{E: primer paso} and Lemma \ref{L: justificaci—n de las cuentas}. Assume now, that $a_1,\dots ,a_k$ have been constructed and set $c_j=a_j^{-\frac 12}$. We have,
\begin{equation*}
F(a_j)=2c_j^{-1}\left(16\, c_j^{-3}-\ddot c_j\right),\ \text{for} \ j=1,\dots k.
\end{equation*}
Let $b_k:[-1,1]\longrightarrow\R$ be the solution to
\begin{equation}\label{E: definicioon}
\begin{cases}
\ddot b_k= -2c_k\left(16\, c_k^{-3}-\ddot c_k\right),\ \text{in}\ [-1,1],\\
b_k(-1)=b_k(1)=0.
\end{cases}
\end{equation} 
Observe that $b_k$ is even, 
\begin{equation}\label{E: formuklaparab}
b_k(t)=2\int_t^1\int_0^s c_k(\tau)\left(16\, c_k^{-3}(\tau)-\ddot c_k(\tau)\right)d\tau ds,\ \text{in}\ [0,1],
\end{equation}
and $\dot b_k<0$ in $(0,1]$. Apply now \eqref{E: convexiad} in  Lemma \ref{L: justificaci—n de las cuentas}, with $a=a_k$ and $b=b_k\xi$, $\xi\in\Rn$. We get 
\begin{equation}\label{E: convexiad2}
\|e^{a_k(t)|x+b_k(t)\xi|^2}u(t)\|\le e^{T_k(t)|\xi|^2+2\|V\|_\infty+\frac 14\|V\|_\infty^2}\sup_{[-1,1]}\|e^{\mu|x|^2}u(t)\|,\ -1\le t\le 1,
\end{equation}
with
\begin{equation*}
\begin{cases}
\partial_t\left(\tfrac 1{a_k}\,\partial_tT_k\right)=-\frac{|\ddot b_k|^2}{F(a_k)},\ \text{in}\ [-1,1],\\
T_k(-1)=T_k(1)=0.
\end{cases}
\end{equation*}
Because $T_k$ is even, the monotonicity of $a_k$, \eqref{E: definicioon} and \eqref{E: formuklaparab}, we get
\begin{equation*}
T_k(t)=2\int_t^1\int_0^s\frac{a_k(s)}{a_k(\tau)}\, c_k(\tau)\left(16\, c_k^{-3}(\tau)-\ddot c_k(\tau)\right)d\tau ds\le  b_k(t),
\end{equation*}
in $[-1,1]$, and from \eqref{E: convexiad2}
\begin{multline}\label{E: bastantebien}
\int_{\Rn}e^{2a_k(t)|x|^2-2|\xi|^2b_k(t)\left(1-a_k(t)b_k(t)\right)+4a_k(t)b_k(t)x\cdot\xi}|u(t)|^2\,dx\\
\le e^{4\|V\|_\infty+\frac 12\|V\|_\infty^2}
\sup_{[-1,1]}\|e^{\mu|x|^2}u(t)\|^2, \ \text{when}\ -1\le t\le 1.
\end{multline}

The latter implies that $u\equiv 0$, when $1-a_k(0)b_k(0)\le 0$, a case in which the process stops. Otherwise, the monotonicity of $a_k$ and $b_k$ implies that $1-a_kb_k>0$, in $[-1,1]$. Multiply then \eqref{E: bastantebien} by $e^{-2\e b_k(t)|\xi|^2}$, $\e>0$, and integrate the corresponding inequality with respect to $\xi$ in $\Rn$. It gives,
\begin{equation}\label{E: bastantebienmejor}
\sup_{[-1,1]}\|e^{a_{k+1}^\e(t)|x|^2}u(t)\|\le \left(1+\tfrac{1}{\e}\right)^{\frac n4}e^{2\|V\|_\infty+\frac 14\|V\|_\infty^2}\sup_{[-1,1]}\|e^{\mu|x|^2}u(t)\|,
\end{equation}
with
\begin{equation*}
a_{k+1}^\e =\frac{\left(1+\e\right)a_k}{1+\e-a_kb_k}\, .
\end{equation*}
Set then
\begin{equation}\label{E: siguientepaso}
a_{k+1} =\frac{a_k}{1-a_kb_k},\quad c_{k+1}=a_{k+1}^{-\frac 12}.
\end{equation}
We get that $a_{k+1}$ is even, $a_{k+1}(1)=\mu$, $a_k<a_{k+1}$, in $(-1,1)$, $\dot a_{k+1}\le 0$, in $[0,1]$
and $F(a_{k+1})>0$, in $[-1,1]$. To verify the latter, recall that
\begin{equation*}
F(a_{k+1})=2c_{k+1}^{-1}\left(16\, c_{k+1}^{-3}-\ddot c_{k+1}\right).
\end{equation*}
From \eqref{E: siguientepaso} and \eqref{E: definicioon}, $c_{k+1}=\left(c_k^2-b_k\right)^{\frac 12}$ and
\begin{equation}\label{E: las derivadas segundas}
\ddot c_{k+1}=c_{k+1}^{-3}\left(16-\tfrac{\dot b_k^2}4+c_k\dot c_k\dot b_k-\dot c_k^2b_k-16\, c_k^{-2}b_k\right).
\end{equation}
Moreover, from \eqref{E: algunas propiedades mas} and \eqref{E: formuklaparab}, $\dot c_k\dot b_k\le 0$ and $16\, b_kc_k^{-2}+\dot b_k^2 >0$, in $[-1,1]$. Thus,
\begin{equation*}
 \ddot c_{k+1}<16\,c_{k+1}^{-3},\ \text{in}\ [-1,1].
\end{equation*}
Also, \eqref{E: bastantebienmejor} implies that
\begin{equation*}
\sup_{[-1,1]}\|e^{\left(a_{k+1}(t)-\e\right)|x|^2}u(t)\|<+\infty,\ \text{for all }\ \e>0,
\end{equation*}
and Lemma \ref{L: justificaci—n de las cuentas} now shows that \eqref{E: algunas propiedades mas}, \eqref{E: algoagradable222} and \eqref{E: algoagradable2222} hold up to $j=k+1$.

When the process is infinite, we have \eqref{E: algoagradable222} and \eqref{E: algoagradable2222} for all $j\ge 1$ and there are two possibilities: either $\lim_{k\to +\infty}a_k(0)=+\infty$ or $\lim_{k\to +\infty}a_k(0)<+\infty$. The first case and \eqref{E: algoagradable222} implies,  $u\equiv 0$, while in the second, the sequence $a_k$ verifies,
\begin{equation}\label{E: acotacion}
\mu\equiv a_1\le a_2\le\dots a_k\le \dots \le \lim_{k\to +\infty}a_k(0),\ \text{in}\ [-1,1], 
\end{equation}
and if $a(t)=\lim_{k\to +\infty}a_k(t)$, set $c=a^{-\frac 12}$. From \eqref{E: siguientepaso}, 
\begin{equation*}
b_k=\frac{a_{k+1}-a_k}{a_k\, a_{k+1}}\, ,
\end{equation*}
$\{b_k\}$ is uniformly bounded and $\lim_{k\to +\infty}b_k(t)=0$, in $[-1,1]$. From \eqref{E: algunas propiedades mas} and the evenness of $c_k$, we have, $\ddot c_k\le 16\, c_k^{-3}$, in $[-1,1]$ and
\begin{equation*}
0\le \dot c_k(t)\le 16\int_0^t c_k^{-3}(s)\,ds,\ \text{in}\ [0,1].
\end{equation*}
Thus, $\{\dot c_k\}$ is uniformly bounded in $[-1,1]$. From, \eqref{E: formuklaparab} and \eqref{E: definicioon}
\begin{align}
&\dot b_k(t)=-\int_0^t32\,c_k^{-2}(\tau)+2\,\dot c_k^2(\tau)\,d\tau+\dot{\widehat{c_k^2}}(t),\label{E: control de la derivada primera}\\
&b_k(t)+\frac 1\mu=\int_t^1\int_0^s 32\, c_k^{-2}(\tau)+2\, \dot c_k^2(\tau)\, d\tau ds +c_k^2(t),\label{E: control de la derivada seguna}
\end{align}
in $[0,1]$, and \eqref{E: control de la derivada primera} shows that $\{\dot b_k\}$ is uniformly bounded in $[-1,1]$, while the uniform boundedness of $\{\ddot c_k\}$ follows from \eqref {E: las derivadas segundas}. Letting $k$ tend to infinity in \eqref{E: control de la derivada seguna}, we find that
\begin{equation*}
\frac 1\mu=c^2(t)+\int_t^1\int_0^s32\, c^{-2}(\tau)+2\, \dot c^2(\tau)\, d\tau ds,\ \text{in}\ [0,1],
\end{equation*}
which implies that
\begin{equation*}
\begin{cases}
\ddot c=16\, c^{-3},\ \text{in}\ [-1,1],\\
\dot c(0)=0, c(1)=\frac 1{\sqrt\mu}\, .
\end{cases}
\end{equation*}
In particular, $a$ is even and
\begin{equation*}
\begin{cases}
\ddot a-\frac{3\dot a^2}{2a\,\,}+32a^3=0,\ \text{in}\ [-1,1],\\
a(1)=\mu.
\end{cases}
\end{equation*}
Because 
\begin{equation*}
\frac{R}{4\left(1+R^2t^2\right)}\, ,\quad R\in\R,
\end{equation*}
are all the possible even solutions of this equation, $a$ must be one of them and 
\begin{equation}\label{E: quien es mu}
\mu =\frac{R}{4\left(1+R^2\right)}\, ,
\end{equation}
for some $R>0$. In particular, $u\equiv 0$, when $\mu >1/8$.

The bounds,
\begin{multline*}
\sup_{[-1,1]}\|e^{a(t) |x|^2}u(t)\|\le e^{2\|V\|_\infty+\frac 14\|V\|_\infty^2}\sup_{[-1,1]}\|e^{\mu |x|^2}u(t)\|,\\
\|\sqrt{1-t^2}\,\,\nabla(e^{(a+\frac{i\dot a}{8a})|x|^2}u)\|_{L^2(\Rn\times [-1,1])}
+\eta_{\e}\|\sqrt{1-t^2}\,e^{\left(a(t)-\e\right)|x|^2}\nabla u\|_{L^2(\Rn\times [-1,1])}\\ \le e^{2\|V\|_\infty+\frac 14\|V\|_\infty^2}\sup_{[-1,1]}\|e^{\mu|x|^2}u(t)\|,
\end{multline*}
follow from \eqref{E: algoagradable2222}.  

The application of the constructive process to the free wave $u_R$ in \eqref{E: el enemigo} shows that $\lim_{k\to +\infty}a_k(0)<+\infty$, when $\mu<\frac 18$ and that the final limit $a$ is determined by the smallest root of the equation \eqref{E: quien es mu}. The same holds, when $\mu=\frac 18$, as it follows by applying the process to the counterexample in the proof of Theorem \ref{T: hardytimeindepent2}.

Theorem \ref{T: lamejora} follows  after undoing the changes of variables at the beginning of the proof and from the bounds in \cite[Theorem 3]{ekpv08b}, when $V$ verifies the first condition or \cite[Theorem 5]{ekpv08b}, when $V$ verifies the second condition in Theorem \ref{T: hardytimeindepent}. The relation between the original $u$ and $v$ is the following:
\begin{multline*}
v(x,t)\\=\left(\tfrac{\sqrt{2\alpha\beta T}}{\alpha(1-t)+\beta (t+1)}\right)^{\frac n2}u\left(\tfrac{\sqrt{2\alpha\beta T}\, x}{\alpha(1-t)+\beta (1+t)}\, , \tfrac{\beta(1+ t)T}{\alpha(1-t)+\beta (1+t)}\right)e^{\left(\alpha-\beta\right) |x|^2/4i\left(\alpha(1-t)+\beta (1+t)\right)}.
\end{multline*}
\end{proof}
\begin{proof}[Proof of Theorem \ref{T: hardytimeindepent2}]
Set
\begin{equation*}
u(x,t)=\left(1+it\right)^{-2k-\frac n2}\left(1+|x|^2\right)^{-k}e^{-\frac{\left(1-it\right)}{4\left(1+t^2\right)}\, |x|^2},
\end{equation*}
for some $\ k>\frac n2$. Then, 
\begin{equation*}
\|e^{|x|^2/8}u(\pm 1)\|<+\infty\quad ,\quad \partial_tu=i\left(\triangle u+V(x,t)u\right),
\end{equation*}
in $\R^{n+1}$, with
\begin{equation*}
V(x,t)=\frac 1{1+|x|^2}\left(\frac{2k}{1+it}+2kn-\frac{4k\left(1+k\right)|x|^2}{1+|x|^2}\right),
\end{equation*}
and
\begin{equation*}
|V(x,t)|\lesssim\frac 1{1+|x|^2}\, ,
\end{equation*}
in $\Rn\times [-1,1]$. What remains follows by modifying the above counterexample with the changes of variables in Lemma \ref{L: transformada}.
\end{proof}
\begin{remark}\label{R: una nota}

The above arguments show that the following also holds under the conditions in Theorem \ref{T: lamejora}: given $\e>0$ there is $\eta_\e$ such that
\begin{multline}\label{E: algo mas de regularidad}
\eta_\e \| \sqrt{t(T-t)}e^{\left(a(t)-\e\right)|x|^2}\nabla u\|_{L^2(\Rn\times [0,T])}\\
\le N\left[\|e^{|x|^2/\beta^2}u(0)\|_{L^2(\Rn)}+ \|e^{|x|^2/\alpha^2}u(T)\|_{L^2(\Rn)}\right],
\end{multline}
with $a$ and $N$ as in Theorem \ref{T: lamejora}. We cannot make $\e=0$ in the exponent of \eqref{E: algo mas de regularidad} because at the end of the process $F(a)$ is identically zero in $[-1,1]$ and we loose the control of $\|x\,e^{a(t)|x|^2}u(t)\|$ in \eqref{E: vaya formulita2}.
\end{remark}
\end{section}
\vfill\eject

\end{document}